\documentclass{lmcs}
\pdfoutput=1

\usepackage{lastpage}
\lmcsdoi{19}{4}{21}
\lmcsheading{}{\pageref{LastPage}}{}{}%
{May~28,~2021}{Dec.~08,~2023}{}


\usepackage{hyperref}
\usepackage[utf8]{inputenc}
\usepackage{graphicx}
\usepackage{amssymb}
\usepackage{amsmath}
\usepackage{amsthm}
\usepackage{latexsym}
\usepackage[all]{xy}
\usepackage[baseline]{euflag}

\def\sub{\subseteq }
\def\overlap{\, \between \,}
\def\sqmeet{>\mkern-13.5mu <}
\newcommand{\primconst}[1] {\mbox{\sf #1}}
\newcommand{\inte}{\primconst{\,int\,}}
\newcommand{\cl}{\primconst{\,cl\,}}
\newcommand{\id}{\primconst{\,id\,}}
\newcommand{\Pos}{\mathrm{Pos}}
\newcommand{\im}{\mathrm{Im}}
\newcommand{\fix}{\mathrm{Fix}}

\begin{document}

\title{Overlap algebras as almost discrete locales}
\thanks{This project has received funding from the European Union's Horizon 2020 research and innovation programme under the Marie Sk\l odowska-Curie grant agreement No 731143 \euflag.}

\author[F.~Ciraulo]{Francesco Ciraulo\lmcsorcid{0000-0002-4957-4799}}
\address{University of Padua, Department of Mathematics ``Tullio Levi-Civita'', Via Trieste 63, 35121 Padova, Italy}
\email{ciraulo@math.unipd.it}

\begin{abstract}
Boolean locales are ``almost discrete'', in the sense that a spatial Boolean locale is just a discrete locale (that is, it corresponds to the frame of open subsets of a discrete space, namely the powerset of a set). This basic fact, however, cannot be proven constructively, that is, over intuitionistic logic, as it requires the full law of excluded middle (LEM). In fact, discrete locales are never Boolean constructively, except for the trivial locale. So, what is an almost discrete locale constructively? Our claim is that Sambin's \emph{overlap algebras} have good enough features to deserve to be called that. Namely, they include the class of discrete locales, they arise as smallest strongly dense sublocales (of overt locales), and hence they coincide with the Boolean locales if LEM holds.
\end{abstract}

\maketitle

\section{Introduction}

A Boolean space is a topological space whose open sets form a complete Boolean algebra (instead of a mere complete Heyting algebra). A Boolean space is almost discrete: a discrete space is just a $T0$ (Kolmogorov) Boolean space or, equivalently, a sober Boolean space. In other words, a spatial Boolean locale \cite{stone_spaces} is the same thing as a discrete locale. For this reason, a Boolean locale can be considered as an ``almost discrete'' locale. All this is fine, at least when classical logic is assumed. Intuitionistically, on the contrary, the property of being Boolean is no longer a necessary condition for discreteness. Indeed, 
if one wants to generalise the above discussion from the usual mathematical universe of sets to the internal universe provided by a topos $\mathcal{E}$, then discrete locales (except for the trivial locale) fail to be Boolean, in general. For instance, the subobject classifier $\Omega$ in $\mathcal{E}$ is a discrete locale (it is the powerset of the terminal object), but it is Boolean precisely when the Law of Excluded Middle (LEM) holds in the internal language of $\mathcal{E}$, that is, when $\mathcal{E}$ is a Boolean topos. Therefore, in a sense, the notions of Boolean-ness and discreteness appear orthogonal to each other from a constructive point of view. So what could serve as an intuitionistic definition of an almost discrete locale? (And so, in a sense, what is a more intuitionistically robust
analogue for Boolean-ness?) We would like our definition of an almost discrete locale to satisfy at least the following two conditions: (i) every discrete locale has to be almost discrete, and (ii) almost discrete locales must coincide with Boolean locales classically. 

The aim of this paper is to show that Sambin's notion of an \emph{overlap algebra} \cite{bp,regolari,formal_regular,cmt,cc} can serve as a constructive definition of an almost discrete locale. Our claim is motivated and supported by some results, the main of which is given in Section~\ref{section o-algebras}: \emph{overlap algebras arise precisely as the smallest strongly dense sublocales \cite{strong_density} of overt locales}. Therefore overlap algebras give an intuitionistically robust account of Boolean-ness. And in the presence of LEM, overlap algebras and Boolean locales coincide (because, in that case, strong denseness boils down to denseness, and overtness comes for free). 

Furthermore, every discrete locale is an overlap algebra, actually an atomic one. This prompts the natural question: what is a ``minimal" condition which makes an overlap algebra discrete? (If we agree that the requiring atomicity seems too strong and not topological at all.) It is natural to wonder if spatiality could be such a property (that is the case in light of LEM). Spatial overlap algebras are sober spaces in which every open set is ``weakly'' regular in the sense that it coincides with the interior of the set of its adherent points. Spaces (not necessarily sober) enjoying this property are studied in Section~\ref{section_spaces} and they are compared with Boolean ones (the two notions coincide classically). Unfortunately, spatial overlap-algebras cannot be proved to be discrete constructively, as we show in Subsection~\ref{counterexample}. So understanding  what the right condition is (which makes an overlap algebra discrete) is still an open question.

Section~\ref{section preliminaries} contains a few preliminaries about locales and a specific limiting result about Boolean locales within an intuitionistic framework (roughly speaking, the result is that there are no nontrivial, overt, Boolean locales). 

Unless otherwise stated, we work constructively in the sense that we assume neither LEM nor the Axiom of Choice.

\section{Preliminaries about frames and locales}\label{section preliminaries}

A {\bf frame} is a complete lattice (that is, a poset with arbitrary joins and meets) satisfying the infinite distributive law
\begin{equation}\label{eq. distributivity}
x\wedge\bigvee_{i\in I}y_i\ =\ \bigvee_{i\in I}(x\wedge y_i)
\end{equation} 
(half of the equality comes for free, of course). A frame homomorphism is a map which preserves finite meets (hence, in particular, the top element 1) and arbitrary joins (hence, in particular, the bottom element 0). Frames are the same thing as complete Heyting algebras, with $x\to y$ = $\bigvee\{z\ |\ z\wedge x\leq y\}$,  although a frame homomorphism need not preserve implications. We often write $-x$ instead of $x\to 0$, the pseudo-complement of $x$. 

The category $\mathbf{Loc}$ of locales is the opposite of the category of frames. For an arrow $f:X\rightarrow Y$ in $\mathbf{Loc}$, we write $\Omega f:\Omega Y\rightarrow\Omega X$ for the corresponding frame homomorphism. 

$\mathbf{Loc}$ has a terminal object $\mathbf{1}$ whose corresponding frame, which is usually written $\Omega$ instead of $\Omega\mathbf{1}$, is the powerset of the set $1=\{0\}$; $\Omega$ can be interpreted as the set of truth-values, and it is in bijection with the two-element set $2=\{0,1\}$ if and only if LEM holds. 

The following are some features of $\Omega$ that are valid intuitionistically. For $p,q\in\Omega$, one has $p\leq q$ precisely when $p=1$ implies $q=1$. Moreover, $p\neq 1$ if and only if $p=0$ if and only if $-p$ = 1.\footnote{On the contrary, $p\neq 0$ implies $p=1$ for all $p\in\Omega$ precisely when $\Omega$ is Boolean.} Also, for $\{p_i\ |\ i\in I\}\subseteq\Omega$, one has $(\bigvee_{i\in I}p_i)=1$ if and only if $p_i=1$ for some $i\in I$. 
 
Given a locale $X$, the unique arrow $!_X:X\to\mathbf{1}$ in $\mathbf{Loc}$ corresponds to the frame homomorphism defined by $\Omega !_X(p)=\bigvee\{x\in\Omega X\ |\ x=1$ and $p=1\}$, for $p\in\Omega$.
A locale $X$ is {\bf overt} \cite{taylor} (or {\bf open} \cite{joyal-tierney}) if $\Omega !_X$ (seen as a monotone map) has a left adjoint $\Pos_X$, that is,
\begin{equation}\label{eq. overtness}
\Pos_X(x)\leq p\textrm{ if and only if }x\leq\Omega !_X(p)
\end{equation}
(which happens precisely when $\Omega !_X$ preserves all meets). In particular, $\Pos_X$ preserves joins. Classically, every locale is overt, and $\Pos_X(x)= 1$ if and only if $x\neq 0$. Constructively $\Pos_X(x)=1$ can be read as a positive (strong) way to express that $x$ is different from 0. We call $\Pos_X$ the {\bf positivity predicate} of $X$.

Every topological space determines a locale $X$, where $\Omega X$ is the frame of open sets. Locales obtained in this way are called {\bf spatial}. They are always overt and $\Pos_X(x)=1$ means that the open set $x$ is inhabited. A locale $X$ is {\bf discrete} if $\Omega X$ is the powerset of some given set; so every discrete locale is spatial and can be seen as the locale corresponding to a discrete space.\footnote{Every discrete locale is (spatial and hence) overt. In fact, overtness appears as one of the two conditions in the characterization of discrete locales in \cite[chapter V, section 5]{joyal-tierney} (the other condition being that the diagonal $X\to X\times X$ is open).}

Different spaces can happen to determine the same spatial locale, the canonical one (up to homeomorphism) being the {\bf sober} one, that is, the one in which every (inhabited) completely prime filter of opens is the collection of open neighbourhoods of a unique point.

\subsection{Sublocales}\label{subsection:sublocales}

A {\bf closure operator} on a poset is an endofunction $c$ such that the conditions $x\leq cx=ccx$ and $x\leq y$ $\Rightarrow$ $cx\leq cy$ hold identically. We write $\fix(c)$ for the collection of all fixed points of $c$; since $c$ is idempotent, $\fix(c)$ = $\im(c)$, the image of $c$. 

A {\bf nucleus} on a locale $X$ is a closure operator $j$ on $\Omega X$ that, in addition, preserves binary meets. In this case, $\fix(j)$ is a frame where finite meets are computed in $\Omega X$ and joins are given by $j$-closure of joins in $\Omega X$. A notable example of a nucleus is the map $x\mapsto--x$ (double negation nucleus).
The locale $X_j$ corresponding to the frame $\Omega X_j$ = $\fix(j)$ is what is called a {\bf sublocale} of $X$ \cite{stone_spaces} and the mapping $x\mapsto jx$ gives a regular monomorphism $X_j\hookrightarrow X$ in $\mathbf{Loc}$.\footnote{The map $x\mapsto jx$ from $\Omega X$ to $\Omega X_j$ is a regular epimorphism of frames: it is the coequalizer of its kernel pair, that is, of the two projections from $\{(x_1,x_2)\in\Omega X\times\Omega X$ $|$ $jx_1=jx_2\}$ (with pointwise operations) to $\Omega X$.} 

For $j_1$ and $j_2$ nuclei on $X$, $X_{j_1}$ is a sublocale of $X_{j_2}$ when $\fix(j_1)\sub \fix(j_2)$ or, equivalently, when $j_2\leq j_1$ pointwise.\footnote{As a sublocale of $X_{j_2}$, $X_{j_1}$ corresponds to the nucleus $j_2x\mapsto j_1x\,(=j_1j_2x)$.} Therefore, the poset of all sublocales of a given locale $X$ is the opposite of the poset $N(X)$ of nuclei on $X$ (with pointwise order); actually, $N(X)$ is a frame \cite{stone_spaces}.
The join of a family $\{X_{j_i}\}_{i\in I}$ of sublocales corresponds to the pointwise meet of nuclei $x\mapsto\bigwedge_{i\in I}(j_ix)$. Meets of sublocales are better seen from another perspective. The sets of the form $\fix(j)$, for $j$ a nucleus on $X$, are precisely the subsets of $\Omega X$ which are closed under arbitrary meets and which contain $x\rightarrow y$ whenever they contain $y$ (see \cite{stone_spaces,pultr, picado-pultr} for details).\footnote{If $S\sub\Omega X$ is such a subset, then the corresponding nucleus is $x\mapsto\bigwedge\{s\in S\ |\ x\leq s\}$.} Therefore, an arbitrary intersection of sets of that form still has the same form, and hence $\bigcap_{i\in I}\fix(j_i)$ corresponds to a sublocale, which is the meet of the family $\{X_{j_i}\}_{i\in I}$.

\paragraph{The sublocale generated by a family of elements.}

Given any $a\in\Omega X$, the nucleus $x\mapsto (x\to a)\to a$ defines the smallest sublocale of $X$ whose frame contains $a$. More generally, 
given any subset $A\sub\Omega X$, the nucleus
\begin{equation}\label{eq def j_P}
j_A(y)\ =\ \bigwedge_{a\in A}\big((y\rightarrow a)\rightarrow a\big)
\end{equation} gives the smallest sublocale of $X$ whose frame contains $A$.
For $x,y\in\Omega X$ we thus have
$x\leq j_A(y)$ if and only if $(y\rightarrow a)\leq(x\rightarrow a)$ for all $a\in A$. Therefore $j_A(y)$ =
$\bigvee\{x\in\Omega X\ |\ (\forall a\in A)(\;y\rightarrow a\,\leq\,x\rightarrow a\;)\}$.

\paragraph{The nucleus ``generated'' by a closure operator.}
Given a closure operator $c:\Omega X\to\Omega X$, one can consider the sublocale of $X$ generated by $\fix(c)$ in the sense of the previous paragraph. The corresponding nucleus is $j_{\fix(c)}(y)$ = $\bigwedge_{b\in\Omega X}((y\rightarrow cb)\rightarrow cb)$. We claim that this can be written as
\begin{equation}\label{eq def j_P bis}
j_{\fix(c)}(y)\ =\ \bigwedge_{a\in\Omega X}\big(a\rightarrow c(a\wedge y)\big)\,.
\end{equation}
Indeed, for every $a\in\Omega X$, $j_{\fix(c)}(y)$ $\leq$ $(y\to c(a\wedge y))\to c(a\wedge y)$ $\leq$ $a\to c(a\wedge y)$ because $a\leq y\to c(a\wedge y)$. Conversely, for every $b\in\Omega X$, $\bigwedge_{a\in\Omega X}(a\rightarrow c(a\wedge y))$ $\leq$ $(y\to cb)\to c((y\to cb)\wedge y)$ $\leq$ $(y\to cb)\to cb$ because $c((y\to cb)\wedge y)=c(y\wedge cb)\leq ccb=cb$.  

Note that $j_{\fix(c)}$ is the nucleus which best approximates $c$ in the following sense: (i) $j_{\fix(c)}(x)\leq c(x)$ for all $x$ in $\Omega X$ and (ii) if $j$ is another nucleus on $X$ such that $j(x)\leq c(x)$ for all $x\in\Omega X$, then $j(x)\leq j_{\fix(c)}(x)$ for all $x\in\Omega X$. 

\subsection{Boolean locales}\label{section boolean}

A locale $X$ is {\bf Boolean} when $\Omega X$ is a Boolean algebra, that is, $x\vee -x=1$ holds identically or, equivalently, $--x=x$ holds identically. 

All possible examples of Boolean locales are of the form $X_{--}$, the sublocale corresponding to the double negation nucleus on some given locale $X$. In fact, $X$ is Boolean if and only if $X_{--}=X$. Moreover (see \cite{stone_spaces} Exercise II.2.4), a sublocale $X_j\hookrightarrow X$ is Boolean if and only if it is generated by a singleton (in the sense of the previous section), that is, there exists $a\in\Omega X$ such that $jx=(x\to a)\to a$ for all $x$.

Classically, every discrete locale is Boolean. Intuitionistically, on the contrary, discrete locales typically fail to be Boolean: an exception is the (locale whose frame is the) powerset of the empty set, that is, the trivial locale. Here we are going to show a particular limiting result about Boolean locales within a constructive framework. 

\paragraph{Fact:} $\Omega$ is the free suplattice\footnote{The category of suplattices (also known as complete join-semilattices) has complete lattices as objects and join-preserving maps as arrows.} on one generator. In other words, for any complete lattice $X$ and any $a\in X$, there is precisely one join-preserving function $f:\Omega\to X$ such that $f(1)=a$. Indeed, every $p\in\Omega$ can be written as $\bigvee\{q\in\Omega\ |\ q=1=p\}$, and hence the only possible candidate for $f$ is given by $f(p)$ = $\bigvee\{f(q)\ |\ q=1=p\}$ = $\bigvee\{r\in X\ |\ r=a\textrm{ and }p=1\}$.

\begin{prop}\label{proposition overt Boolean}
If there exists a locale $X$ with the following two properties
\begin{enumerate}
\item $X$ is Boolean, and
\item there exists a join-preserving map $F:\Omega X\rightarrow\Omega$ such that $F(1)=1$,
\end{enumerate}
then $\Omega$ is Boolean, that is, LEM holds.
\end{prop}
\begin{proof}
The map $F\circ\Omega !_X:\Omega\to\Omega$ preserves joins and top, and hence it is the identity by the fact recalled above. Moreover, $F(-\Omega !_X(p))$ = $-p$, because $-p$ = 1 iff $p=0$ iff $p\neq 1$. Therefore $p\vee-p$ = $F\big(\Omega !_X(p)\big)\vee F\big(-\Omega !_X(p)\big)$ = $F\big(\Omega !_X(p)\vee-\Omega !_X(p)\big)$ = $F(1)$ = 1 for every $p\in\Omega$.
 \end{proof}

In view of this result, there are a number of things one cannot expect to prove within an intuitionistic setting, in general. For instance, it is impossible to construct a locale $X$ which is Boolean, overt and with $\Pos_X(1)=1$. Also, it is impossible to construct any point of any Boolean locale.\footnote{See \cite[Proposition 4.3]{curi} for a similar limiting result.} 

This limitation appears still more striking within a topos $\mathcal{E}$ in which LEM is provably false: if $\neg(\forall p\in\Omega)(p\vee-p)$ is provable in the internal logic of $\mathcal{E}$,\footnote{This happens, for instance, in the topos of sheaves over a $T1$ space without isolated points.} then the only overt Boolean locale in $\mathcal{E}$ is the trivial locale (because $\Pos_X(1)\neq 1$ means 1=0).

A topological space is {\bf Boolean} if the corresponding locale is Boolean, that is, if its open subsets form a Boolean algebra. It is a corollary of the above discussion that the existence of an inhabited and Boolean topological space is equivalent to LEM.

\section{Overlap algebras}\label{section o-algebras}

Let $\inte$ be the interior operator on the subsets of a topological space $X$. An open set $A$ is {\bf regular} if it equals
the interior of its closure, that is, $A$ $=$ $\inte(\inte A^c)^c$, where $(\ )^c$ is the set-theoretic complement. Since $\inte A^c$ is just the pseudo-complement $-A$ in the frame of open sets, we have that $A$ is regular when $A=--A$. For this reason, an element $x$ of a frame/locale is said {\bf regular} if $x=--x$; and a locale is Boolean if and only if all its elements are regular.

Given a set $D$ of points, there are at least two (classically equivalent) ways to define the topological closure of $D$. One possibility is to take $(\inte D^c)^c$. The other possibility is to consider the set $\cl D$ of all adherent points of $D$; so $x\in\cl D$ if and only if $x\in A\Rightarrow D\overlap A$ for every open $A$. Here, following Sambin, we write $X\overlap Y$ to mean that $X\cap Y$ is inhabited (classically, $X\cap Y\neq\emptyset$), that is, $X$ and $Y$ \emph{overlap} each other. 
The latter definition, which is the one commonly adopted in a constructive approach, results in an intuitionistically weaker notion of closure\footnote{There exists an analogous notion of weak closure for sublocales \cite{strong_density} that is defined in terms of the notion of strong density recalled in Subsection~\ref{section density} below. For our purposes, however, we do not need to recall explicitly such a notion of weak closure for sublocales (thankfully, since an intrinsic characterization of weakly closed nuclei is not at hand \cite[p.7]{strong_density}).} in the sense that $D\sub \cl D\sub (\inte D^c)^c$ always holds, and hence  $D=(\inte D^c)^c$ implies $D=\cl D$. Therefore, every closed subset is weakly closed; and the converse holds only classically.\footnote{Indeed, $\cl D$ = $D$ = $\inte D$ holds for all $D$ in a discrete space (such as $\Omega$); therefore the statement ``$(\inte D^c)^c=\cl D$ holds identically in all spaces'' is precisely $LEM$.} 

By replacing $(\inte A^c)^c$ with $\cl A$ we get an intuitionistically weaker notion of regularity for open sets. We shall use the term {\bf weakly regular} for an open set $A$ such that $A=\inte\cl A$. 

There are spaces (notably the discrete ones) in which all open subsets are weakly regular, so that the following equation (between operators on subsets) holds:
\begin{equation}\label{eq._ici=i}
\inte\cl\inte=\inte\,.
\end{equation} A space of this kind is a ``positive'' version of a Boolean space (perhaps the name ``weakly Boolean" would be appropriate); it corresponds, as we will see, to a spatial overlap algebra.

Here we show how to express ``weakly regular'' in the point-free language of overt locales. As usual, we seek inspiration in the spatial case.\\ 

For $X$ a spatial locale and $a\in\Omega X$, the open set $\inte\cl a$ is the union of all opens $x\in\Omega X$ such that $x\sub\cl a$. By definition, $x\sub\cl a$ means that $z\overlap x$
implies $z\overlap a$ for all $z\in\Omega X$. Now $x\overlap y$ holds precisely when $\Pos_X(x\wedge y)=1$ (spatial locales are overt). Thus, we propose the following.

\begin{defi}\label{def weak regular}
Let $X$ be an overt locale with positivity predicate
$\Pos_X$. We say that $a\in\Omega X$ is {\bf weakly regular} when
\begin{equation}\label{eq. def. weakly regular}
a\ =\ \bigvee\big\{x\in\Omega X\ |\ (\forall z\in\Omega X)\big(\Pos_X(z\wedge x)\leq\Pos_X(z\wedge a)\big)\big\}\ .
\end{equation} 
In other words, $a\in\Omega X$ is {\bf weakly regular} when, for every $x\in\Omega X$,
\begin{equation}
if\ \Pos_X(z\wedge x)\leq\Pos_X(z\wedge a)\ for\ all\ z\in\Omega X,\ then\ x\leq a\ .
\end{equation}
\end{defi}

We shall show  (Proposition~\ref{thm main}) that the weakly regular elements of an overt locale $X$ form a sublocale which is, in fact, the smallest strongly dense sublocale of $X$.\footnote{The notion of strong density (and a corresponding notion of weak closedness) where introduced in \cite{strong_density}; it is recalled in Subsection~\ref{section density} where it plays a fundamental role.} 

\begin{defi}\label{def o-algebra}
An {\bf overlap algebra} ({\bf o-algebra} for short) is an overt locale in which all elements are weakly regular.
More explicitly, an o-algebra is an overt locale $X$ such that, for every $x,y\in\Omega X$,
\begin{equation}\label{eq. def. o-algebra}
if\ \Pos_X(z\wedge x)\leq\Pos_X(z\wedge y)\ for\ all\ z\in\Omega X,\ then\ x\leq y\ .
\end{equation}
\end{defi}

The following Proposition is evidence of the fact that an o-algebra is an appropriate intuitionistic version of a Boolean locale, that is, of an almost discrete locale. Further evidence will be provided by Corollary~\ref{corollary_o-algebra}, which we consider the main result of this paper.

\begin{prop}\label{prop o-algebra vs boolean} Let $X$ be a locale.
\begin{enumerate}
\item If $X$ is discrete, then  $X$ is an o-algebra.\footnote{Actually, every overt sublocale of a discrete locale is an o-algebra~\cite[Corolalry 5.7]{cc}. In general, we refer the reader to~\cite[Section 5.7]{cc} for more on sublocales of an o-algebra. There, in particular, it is shown that every open sublocale of an o-algebra is an o-algebra as well.}
\item If $X$ is overt and Boolean, then $X$ is an o-algebra.
\item Classically, if $X$ is an o-algebra, then $X$ is (overt and) Boolean.
\end{enumerate}
\end{prop}
\begin{proof}
$1.$ Let $X$ be the powerset of $S$. In this case $\Pos_X(x)=1$ means that $x$ is an inhabited subset of $S$. Thus, $\Pos_X(x\cap y)=1$ is just $x\overlap y$ and condition \eqref{eq. def. o-algebra} reads
$\forall z\big((z\overlap x)\Rightarrow (z\overlap y)\big)\Longrightarrow(x\sub y)$, which is true (make $z$ vary over all singletons).\\
$2.$ Assume $\forall z(\Pos_X(z\wedge x)\leq\Pos_X(z\wedge y))$; in particular, $\Pos_X(-y\wedge x)\leq\Pos_X(-y\wedge y)$ = $\Pos(0)=0$, that is, $-y\wedge x\leq\Omega !(0)=0$. Hence $x\leq--y$ and so $x\leq y$ as wished.\\
$3.$ Classically, $\Pos_X(x)=1$ is $x\neq 0$. Thus the antecedent of condition \eqref{eq. def. o-algebra} becomes $\forall z((z\wedge x\neq 0)\Rightarrow(z\wedge y\neq 0))$ which is classically equivalent to $\forall z((z\wedge y=0)\Rightarrow(z\wedge x=0))$. This means precisely $\forall z((z\leq-y)\Rightarrow(z\leq -x))$, that is, $-y\leq-x$. So \eqref{eq. def. o-algebra} states that $(-y\leq-x)\Rightarrow(x\leq y)$ for all $x,y\in\Omega X$. Therefore, $X$ is Boolean.
 \end{proof}

Item $1.$ says that every powerset is an o-algebra in a natural way: this is the prototypical example which actually motivated the introduction of o-algebras by Sambin (see \cite{bp} and also \cite{regolari,formal_regular,cmt}). Items $2.$ and $3.$ together say that o-algebras and Boolean locales coincide classically; so spatial o-algebras correspond to discrete spaces, classically, as a consequence of Proposition~\ref{prop almost discrete} below. Item $3.$ cannot hold intuitionistically because discrete locales are o-algebras by item $1.$, but they are not Boolean unless LEM holds. Item $2.$ is of questionable interest, of course, because of Proposition~\ref{proposition overt Boolean}; indeed, if LEM fails, then the only overt and Boolean locale is the trivial one (that is, the powerset of the empty set).\\

O-algebras are often (and originally were) presented as complete lattices equipped with an extra binary relation $\sqmeet$ (intended as an algebraic version of the relation $\overlap$ between subsets), essentially as follows. 

\begin{prop}\label{prop. def. o-algebra}
The frame underlying an o-algebra is precisely a complete lattice equipped with a symmetric relation $\sqmeet$ such that the following conditions are identically satisfied.
\begin{subequations}
\begin{align}
(x\wedge z)\sqmeet y & \qquad\textrm{if and only if}\qquad x\sqmeet(z\wedge y) \label{eq. transfer} \\
x\sqmeet\left(\bigvee_{i\in I}y_i\right) & \qquad\textrm{if and only if}\qquad (\exists i\in I)(x\sqmeet y_i) \label{eq. splitting} \\
x\leq y & \qquad\textrm{if and only if}\qquad \forall z(z\sqmeet x\Rightarrow z\sqmeet y) \label{eq. density}
\end{align}
\end{subequations}
\end{prop}
\begin{proof} Assume we have a complete lattice satisfying the three conditions above. First, we show that such a lattice satisfies the infinite distributive law \eqref{eq. distributivity} so that it is a frame. By \eqref{eq. density}, that is equivalent to show that $z\sqmeet(x\wedge\bigvee_{i\in I}y_i)$ $\Rightarrow$ $z\sqmeet(\bigvee_{i\in I}(x\wedge y_i))$ for all $z$. So assume the premise; by \eqref{eq. transfer} one gets $(z\wedge x)\sqmeet(\bigvee_{i\in I}y_i)$ and hence $(z\wedge x)\sqmeet y_i$ for some $i\in I$, by \eqref{eq. splitting}. By \eqref{eq. transfer}, that becomes $z\sqmeet(x\wedge y_i)$ for some $i\in I$ and so $z\sqmeet(\bigvee_{i\in I}(x\wedge y_i))$ by \eqref{eq. splitting}. \\
Second, we show that the corresponding locale $X$ is overt with $\Pos_X(x)=1$ $\Leftrightarrow$ $x\sqmeet x$. We must check that $x\sqmeet x$ $\Rightarrow$ $p=1$ if and only if $x\leq\Omega !_X(p)$. By \eqref{eq. density} and \eqref{eq. splitting}, $x\leq\Omega !_X(p)$ means that, for every $z$, if $z\sqmeet x$, then both $z\sqmeet 1$ and $p=1$. Clearly, $z\sqmeet x$ yields $z\sqmeet 1$; therefore $x\leq\Omega !_X(p)$ is just equivalent to $\forall z(z\sqmeet x\Rightarrow p=1)$. By logic, this is just $\exists z (z\sqmeet x)$ $\Rightarrow$ $p=1$. Now $\exists z (z\sqmeet x)$ is tantamount to $x\sqmeet x$ and we are done.\\
Third, we show that \eqref{eq. def. o-algebra} holds for $\Pos_X$. By \eqref{eq. transfer}, $(x\wedge y)\sqmeet(x\wedge y)$ is equivalent to $x\sqmeet y$. So $\Pos_X(x\wedge y)=1$ is equivalent to $x\sqmeet y$ and hence \eqref{eq. def. o-algebra} follows by \eqref{eq. density}.\\
We now come to the opposite direction. Let $X$ be an o-algebra. We define $x\sqmeet y$ as $\Pos_X(x\wedge y)=1$. Clearly, $\sqmeet$ is symmetric and satisfies \eqref{eq. transfer}. Also, \eqref{eq. splitting} easily follows from \eqref{eq. distributivity} and from the fact that $\Pos_X$ preserves joins. Finally, \eqref{eq. density} is a consequence of \eqref{eq. def. o-algebra} and of the fact that $\Pos_X$ is monotone.
 \end{proof}

\subsection{O-algebras are smallest strongly dense sublocales}\label{section density}

There exists a well-known connection between Boolean locales and dense sublocales. A sublocale $X_j\hookrightarrow X$ is {\bf dense} if $\Omega X_j$ contains the bottom element of $\Omega X$, that is, if $j(0)=0$. And $X_{--}$ is the smallest dense sublocale of $X$. Now $X_{--}$ is Boolean, and every Boolean locale is of this form (we refer the reader to \cite{stone_spaces} for details). Therefore, Boolean locales arise precisely as smallest dense sublocales. Here we want to show that o-algebras enjoy an analogous characterization.

A stronger notion of density was introduced in \cite{strong_density}: $X_j\hookrightarrow X$ is {\bf strongly dense} if $j(\Omega !_X(p))$ = $\Omega !_X(p)$ for all $p\in\Omega$, that is, $\Omega !_{X_j}(p)$ = $\Omega !_X(p)$ for all $p\in\Omega$.
In particular, $j(0)=j(\Omega !_X(0))=\Omega !_X(0)=0$. Therefore, strong density implies density. Classically, the converse holds as well. 

Thus, $X_j$ is strongly dense precisely when $Im(\Omega !_X)\sub \fix(j)$.
Therefore, every locale has a smallest strongly dense sublocale, namely, the sublocale generated by $Im(\Omega !_X)$ as described in Section~\ref{subsection:sublocales}. Thus, according to equation \eqref{eq def j_P}, the nucleus corresponding to the smallest strongly dense sublocale is given by
\begin{equation}\label{eq. def. min. strongly dense}
y\longmapsto\bigwedge_{p\in\Omega}\Big(\big(y\rightarrow\Omega !_X(p)\big)\rightarrow\Omega !_X(p)\Big)\ .
\end{equation}
Locales in which this nucleus is the identity are constructive versions of Boolean locales. It is not hard to show that every discrete locale is of this type. Therefore, the locales in which \eqref{eq. def. min. strongly dense} is the identity could be called ``almost discrete''. 
However, since every discrete locale has to be overt, we prefer to include overtness among the properties that an almost discrete locale should enjoy. 

So, from now on, we switch our attention to strongly dense sublocales that are also overt, and we look for an overt version of \eqref{eq. def. min. strongly dense}. Let us start with the following fact \cite[Lemma 1.11]{strong_density}.

\begin{lem}\label{lemma X_j weakly dense -> X dense} Let $X_j\hookrightarrow X$ be a sublocale.
\begin{enumerate}
\item If $X_j$ is strongly dense, then $X$ is overt precisely when $X_j$ is overt.
\item If $X$ is overt, then $X_j$ is strongly dense precisely when $\Pos_X\circ j=\Pos_X$.\footnote{So, in the category of overt locales, strong density becomes equivalent to $\Pos_X\circ j = \Pos_X$. This is a special case of the definition of density proposed in \cite{townsend} for an abstract category resembling $\mathbf{Loc}$.}
\end{enumerate}
\end{lem}
\begin{proof} $1.$
Let $X$ be overt. For $jx\in\Omega X_j$ and $p\in\Omega$ one has: $jx\leq\Omega !_{X_j}(p)=\Omega !_X(p)$ iff $\Pos_X(jx)\leq p$. Therefore, $\Pos_{X_j}$ exists and is the restriction of $\Pos_X$ to $\Omega X_j$. Conversely, let $X_j$ be overt. For $x\in\Omega X$ and $p\in\Omega$ one has: $\Pos_{X_j}(j(x))\leq p$ iff $j(x)\leq \Omega !_{X_j}(p)=j(\Omega !_X(p))$ iff $x\leq j(\Omega !_X(p))=\Omega !_X(p)$. Hence, $X$ is overt with $\Pos_X(x)=\Pos_{X_j}(jx)$.

$2.$ Since $X$ is overt, $x\leq\Omega !_X(\Pos_X(x))$ for every $x\in\Omega X$. Now if $X_j$ is strongly dense, then $jx\leq j(\Omega !_X(\Pos_X(x)))$ = $\Omega !_X(\Pos_X(x))$, and hence
$\Pos_X(jx)\leq\Pos_X(x)$, that is, $\Pos_X\circ j =\Pos_X$.
Conversely, given $p\in\Omega$, we must check that $j(\Omega!_X(p))\leq\Omega !_X(p)$, that is, $\Pos_X(j(\Omega !_X(p)))\leq p$, which follows from $\Pos_X\circ j$ = $\Pos_X$ and overtness of $X$. 
 \end{proof}

In particular, the following facts are equivalent for a locale $X$:
\begin{enumerate}
\item $X$ is overt;
\item every strongly dense sublocale of $X$ is overt;
\item there exists an overt, strongly dense sublocale of $X$.
\end{enumerate}
And the smallest overt, strongly dense sublocale of $X$ exists if and only if $X$ is overt and, in that case, it is just the smallest strongly dense sublocale of $X$. The following gives an alternative formula for the nucleus defined in \eqref{eq. def. min. strongly dense} in the case of an overt locale.

\begin{prop}\label{thm main}
For $X$ an overt locale, the smallest overt, strongly dense sublocale of $X$ corresponds to the nucleus
\begin{equation}\label{eq. def. j_exists}
R_X(y)\ =\ \bigvee\big\{x\in\Omega x\ |\ 
(\forall z\in\Omega X)\big(\Pos_X(z\wedge x)\leq\Pos_X(z\wedge y)\big)
\big\}\ .
\end{equation}
\end{prop}
\begin{proof}
Let $X$ be an overt locale. By triangular identities for adjunctions, $\Omega !_X$ = $\Omega !_X\circ\Pos_X\circ\Omega !_X$ and hence $Im(\Omega !_X)$ = $\fix(\Omega !_X\circ\Pos_X)$. Therefore, the smallest strongly dense sublocale of $X$, which is the sublocale generated by $Im(\Omega !_X)$, coincides with the sublocale generated by $\fix(\Omega !_X\circ\Pos_X)$. By \eqref{eq def j_P bis}, the corresponding nucleus is
$
y\mapsto\bigwedge_{z\in\Omega X}\big(z\rightarrow\Omega !_X\circ\Pos_X(z\wedge y)\big)$.
We call it $R_X$. Thus
$x\leq R_X(y)$ if and only if $\Pos_X(z\wedge x)\leq\Pos_X(z\wedge y)$ for all
$z\in\Omega X$. Thereby \eqref{eq. def. j_exists} follows.
 \end{proof}

In view of Definitions~\ref{def weak regular} and \ref{def o-algebra}, Proposition~\ref{thm main} says that the weakly regular elements of $X$ are just the fixed points of $R_X$, and that o-algebras are precisely the overt locales $X$ for which the nucleus $R_X$ is the identity. Thus, we obtain the following at once.
\begin{cor}\label{corollary_o-algebra} The smallest (overt) strongly dense sublocale of an overt locale is an o-algebra and every o-algebra can be obtained in this way.
\end{cor}
\begin{cor}\label{corollary characterization o-algebras} An overt locale $X$ is an o-algebra if and only if it coincides with its own smallest (overt) strongly dense sublocale. 
\end{cor}

Over classical logic, as we have already noted, o-algebras boil down to Boolean locales, and $R_X$ becomes just the double negation nucleus.
Intuitionistically, one only has 
 $R_X(y)\leq --y$ for all $y$ and hence the double-negation sublocale $X_{--}$ is in fact a sublocale of $X_{R_X}$. 

\section{Discrete spaces and spatial o-algebras}\label{section_spaces}

Aim of this section is to investigate the connection between spatial o-algebras and discrete spaces. 

Spatial o-algebras can be identified with those sober spaces that satisfy the identity $\inte\cl\inte=\inte$ (all open subsets are weakly regular). Classically, such spaces are Boolean and $T0$, hence discrete (see Proposition~\ref{prop almost discrete} below). So the question arise whether spatial o-algebras are discrete or not intuitionistically. Unfortunately, the answer is negative in view of the Brouwerian counterexample provided below. Thus the problem remains of finding an additional requirement on a spatial o-algebra that would make it discrete.\footnote{It is known (see \cite[Proposition 1.2]{cc}, for instance) that every atomic locale is discrete (and vice versa). In particular, every atomic (and spatial) o-algebra is discrete. Requiring atomicity, however, is not so satisfactory for our purpose: it seems too strong (if compared to the classical case, in which spatiality is enough) and not really topological in nature.}
We start with a constructive characterization of discrete spaces. 

A topological space is Boolean if the Heyting pseudo-complement on open sets $\inte(\_)^c$ is an involution. Classically, this is equivalent to saying that every closed set is open. Accordingly, we can consider  two corresponding constructive conditions, which are obtained by writing $\cl$ in place of $(\inte(\_)^c)^c$, namely,
$\inte\cl\inte=\inte$ and $\inte\cl=\cl$.

In what follows, a space $X$ is said to be 
\begin{description}
\item[$R0$] if $y\in\cl\{x\}$ $\Longrightarrow$ $x\in\cl\{y\}$ for all $x,y\in X$;
\item[$T0$] if $\cl\{x\}=\cl\{y\}$ $\Longrightarrow$ $x=y$ for all $x,y\in X$; 
\item[$T1$] if $\cl\{x\}=\{x\}$ for all $x\in X$.
\end{description}
Note that a space which is both $T0$ and $R0$ is also $T1$ (and vice versa).\footnote{Indeed, from $y\in\cl\{x\}$ one gets both
$\cl\{y\}\sub\cl\{x\}$ and $\cl\{x\}\sub\cl\{y\}$, and hence $y=x$; so $y\in\cl\{x\}$ implies $y\in\{x\}$ for every $y$, that is, $\cl\{x\}\subseteq\{x\}$ for all $x$.}

\begin{prop}\label{prop almost discrete}
For a topological space $X$ the following are equivalent:
\begin{enumerate}
\item $X$ is discrete (that is, $\inte=\id$);
\item $X$ is $T0$ and satisfies $\inte\cl=\cl$.
\end{enumerate}
\end{prop}
\begin{proof}
The proof of $2\Rightarrow 1$ follows by the following three facts together.
\begin{itemize}
\item If $\inte\cl=\cl$, then $\cl\inte=\inte$.\\
We must show that $\cl A\sub A$ for every open set $A$. Thus let $x\in\cl A$. The set $\cl\{x\}$ is open, and hence it is an open neighbourhood of $x$ by assumption; therefore, it must be the case that $\cl\{x\}\overlap A$, which implies $x\in A$.
\item If $\cl\inte=\inte$ and $X$ is $T0$, then $X$ is $R0$ (and hence $T1$).\\ 
Let $y\in\cl\{x\}$. For every open $A\ni x$, we have $y\in\cl\{x\}\sub\cl A=A$. Thus, $\{y\}\overlap A$. Therefore, $x\in\cl\{y\}$.
\item If $\inte\cl=\cl$ and $X$ is $T1$, then $X$ is discrete.\\ 
For every $x\in X$ we have $\inte\{x\}$ = $\inte\cl\{x\}$ = $\cl\{x\}$ = $\{x\}$. So every subset of $X$ is open.
\end{itemize}
The proof of $1\Rightarrow 2$ is trivial, of course.\qedhere
\end{proof}

\subsection{A Brouwerian counterexample}\label{counterexample}

Thanks to the proof of Proposition~\ref{prop almost discrete}, the following implication holds.
\begin{equation}\label{eq. forms of almost discreteness}
\inte\cl=\cl\quad\Rightarrow\quad\inte\cl\inte=\inte
\end{equation}
Classically, the converse implication holds as well, and any of them can be used to define a Boolean space. 
Classically, therefore, a discrete space is just a Kolmogorov, Boolean space or, equivalently, a sober Boolean space (that is, a spatial Boolean locale). 

Constructively, the picture is quite different, as we now see.

\begin{lem}
For each $p\in\Omega$ there is a topological space $(2,\tau_p)$ on the set $2=\{0,1\}$ such that:
\begin{enumerate}
    \item every subspace of $(2,\tau_p)$ is weakly closed (that is, $\cl=\id$);
    \item $(2,\tau_p)$ is discrete if and only if  $p\vee\neg p=1$.
\end{enumerate}
\end{lem}
\begin{proof}
Given any subset $P\sub 2$, the family $$B_P=\big\{\{0\}\cap P,\{1\}\cap P,\{0\}\cup P,\{1\}\cup P\big\}$$ is a base for a topology on 2. Indeed, 2 = $(\{0\}\cup P)\cup(\{1\}\cup P)$; and the intersection of any two elements of $B_P$ is either empty or an element of $B_P$ or $P$, which is the union of the two elements $\{0\}\cap P$ and $\{1\}\cap P$ of $B_P$.\footnote{All subsets and supersets of $P$ are open in such a topology; classically, there is no other open set.}
Note that if $P=2$ or $P=\emptyset$, then the corresponding topology is discrete. Conversely, if the topology is discrete, then $\{0\}$ would be open, that is, there should exist $A\in B_P$ such that $0\in A\sub\{0\}$; so either $0\in P$ or $P\subseteq\{0\}$. Similarly, either $1\in P$ or $P\subseteq\{1\}$, because $\{1\}$ would be open. Therefore, either $P=2$ or $P=\emptyset$.

From now on, let $P$ be $\{x\in 2\ |\ p=1\}$ and consider the topology $\tau_p$ generated by the base 
$B_P$. 

By the previous discussion, $\tau_p$ is discrete if and only if $p\vee\neg p=1$.

We claim that every set $D$ is weakly closed in $(2,\tau_p)$ whatever $p$ is. Let $x\in\cl D$; we must show that $x\in D$. Choose $E=\{x\}\cup P$, which is an open neighbourhood of $x$; therefore $E\overlap D$. Thus either $x\in D$ or $P\overlap D$. In the former case we are done. In the latter case, $p$ holds and so $\{x\}$ = $\{x\}\cap P$ is another open neighbourhood of $x$; therefore $x\in D$.
 \end{proof}

\begin{prop}
The following are equivalent:
\begin{enumerate}
    \item LEM;
    \item for all topological spaces, $\cl=\id$ $\Rightarrow$ $\inte=\id$  
(if every subset is weakly closed, then the space is discrete);
    \item for all topological spaces, $\inte\cl\inte=\inte$ $\Rightarrow$ $\inte\cl=\cl$ (if every open set is weakly regular, then every weakly closed set is open);
    \item every spatial o-algebra is discrete.
\end{enumerate}
\end{prop}
\begin{proof}
Consider the family of topological spaces $(2,\tau_p)$ with $p\in\Omega$ as in the previous lemma. As $\cl=\id$ implies $\inte\cl\inte=\inte$, every $\tau_p$ is a spatial o-algebra. Thus, the equivalence between the above items follows easily from the previous Lemma. (Show that $4.\Rightarrow 3.\Rightarrow 2.\Rightarrow 1.$ and recall that LEM makes o-algebras and Boolean locales coincide.) 
\end{proof}

\section*{Conclusions}

The class of Boolean locales has some disadvantages intutionistically, for instance, it does not contain discrete spaces. On the other hand, the class of ``overlap algebras'' (a particular kind of overt locales) has some interesting features: (i) every discrete locale is an overlap algebra; (ii) the smallest strongly dense sublocale of an overt locale is an overlap algebra, and every overlap algebra arises in this way; (iii) under classical logic, an overlap algebra is the same as a Boolean locale. 

Spatial overlap algebras correspond to sober spaces in which every open is ``weakly'' regular, that is, the equation $\inte\cl\inte=\inte$ is identically satisfied. Classically, these are precisely the discrete spaces. Thus, spatial o-algebras could deserve the name of ``almost discrete locales''.
 However, spatial overlap algebras cannot be proven to be discrete without assuming the law of excluded middle. 
Finding a ``minimal'' additional condition that would  make a spatial o-algebra discrete remains an open problem.

 \section*{Acknowledgements} The author wishes to thank three anonymous reviewers who, through their kind comments, encouraged him to improve the paper and make it more readable.

\bibliographystyle{alphaurl}
\bibliography{references}

\end{document}